\documentclass[12pt,a4paper]{article}
\usepackage{graphicx}
\usepackage{german}
\usepackage{amsmath}
\usepackage{amssymb}
\usepackage{enumerate}
\usepackage{latexsym}

\begin{document}
	\title{\begin{figure}[htp]
		\centering
		\includegraphics[scale=1.0]{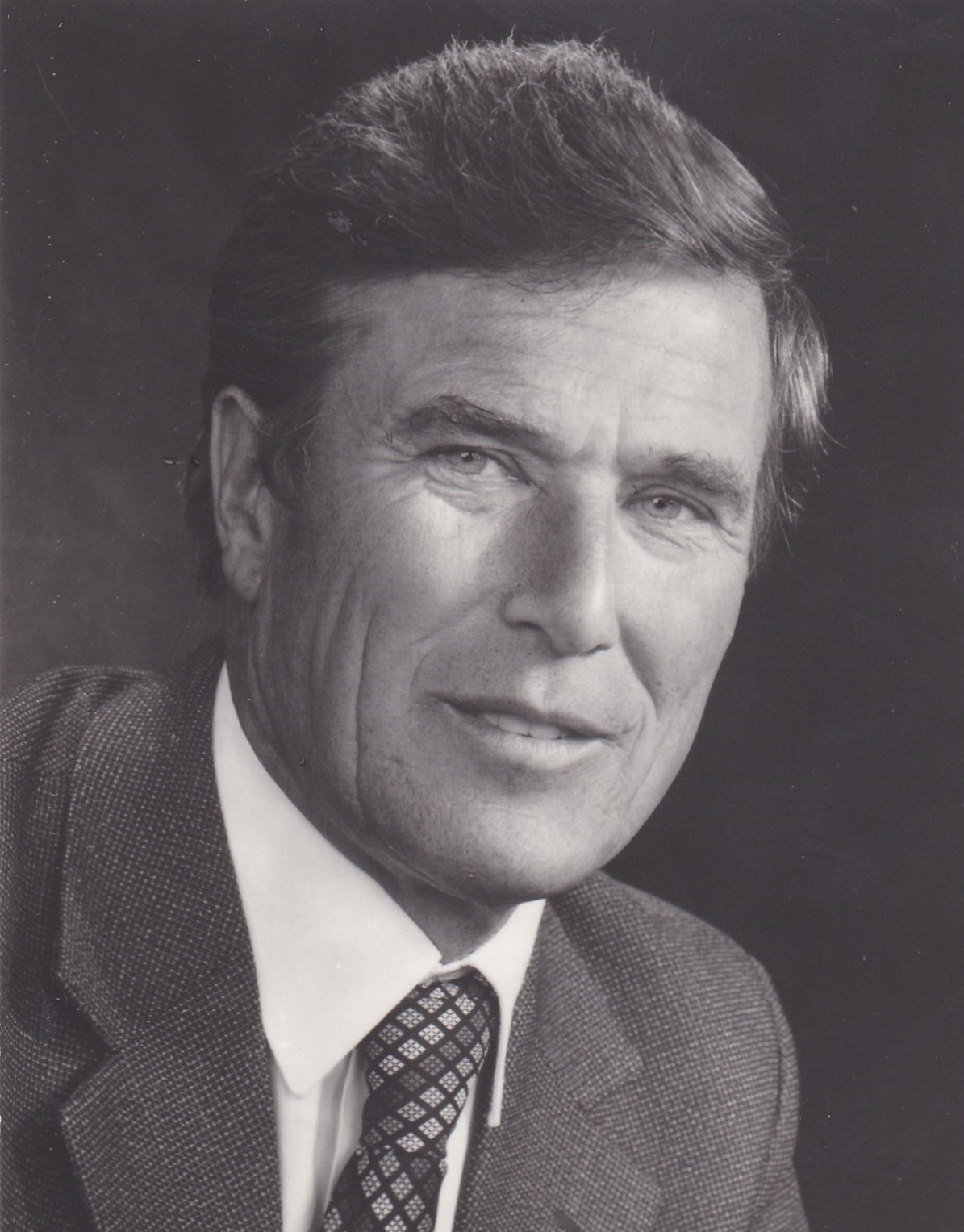}
	\end{figure}
    Helmut Karzel \, (1928--2021)}
\author{Hans Havlicek, Alexander Kreuzer, Hans-Joachim Kroll\\ and Kay S\"{o}rensen}

	\date{} 	
	\maketitle

%%%%%%%%%%%%%%%%%%%%%%%%%%%%%%%%%%%%%%%%%%%%%%%%%%%%%%%%%%%%%%%%%%%%%%%%%%%%%%%%

We are saddened to report that Professor Dr.\,Dr.\,h.c.\ Helmut Karzel passed
away on June 22, 2021, in Burgwedel, Germany, at the age of 93. He will always
be in our best memories as teacher, friend and valued member of the worldwide
community of mathematicians. Our deepest sympathy goes out to his three
children D\"{o}rte, Herbert, Barnim and to Gisela, his partner over the last years.
\par

Helmut Karzel was born on January 15, 1928, in Sch\"{o}neck in Westpreu{\ss}en, Poland.
He spent his youth in Posen, where he attended school from 1934 to 1945. During
the last year of World War II he had to serve as flak helper. At the end of the
war he was advised to go west. So, in 1945 he reached the German town of
Magdeburg, where an uncle lived. There, in 1946, he took the ``Abitur'' at the
Otto von Guericke Schule. His mother and two of his three siblings had to flee
from Posen and all of them met again in Magdeburg. \par

By crossing several occupation zone boundaries (at night), Karzel reached
Freiburg, Germany, in 1947. He received the admission to study mathematics and
physics at the University of Freiburg. Prior to his studies, he was obliged to
help for eight weeks in clearing up the rubble of the bomb attacks. Among
Karzel's teachers in Freiburg was Emanuel Sperner, who became Professor at the
University of Bonn in 1950. Karzel moved to Bonn and obtained his doctorate
under Sperner's supervision in 1951. That year he also received the Hausdorff
Memorial Prize, which is awarded annually by the University of Bonn for the
best dissertation in mathematics of the past academic year. \par

After his studies, Karzel worked as ``Assistent'' at the University of Bonn
and, starting from 1954, at the University of Hamburg. Karzel accomplished his
``Habilitation'' in the year 1956. He stayed as Visiting Associate Professor at
the University of Pittsburgh, Pennsylvania, USA, in 1961/62. Furthermore, he
was on leave from Hamburg for a guest professorship at the Technische
Universit\"{a}t Karlsruhe from 1967 to 1968. During the summer term 1968 he taught
in Hamburg as well as in Karlsruhe. In 1968, Karzel was appointed Chair of
Geometry at the Technische Hochschule Hannover, Germany. A few years later, in
1972, he accepted an offer from Technische Universit\"{a}t M\"{u}nchen as Chair of
Geometry, a position he retained until his retirement in 1996. \par

A major aim of Helmut Karzel's work was to maintain contact with colleagues
worldwide. Consequently, Karzel lectured as Visiting Professor at numerous
universities outside Germany: Bologna, Italy (Adriano Barlotti); Brescia, Italy
(Mario Marchi); Toronto, Canada (Erich Ellers); College Station, Texas, USA
(Carl J.~Max\-son); Rome, Italy (Giuseppe Tallini); Teesside, England (Allan
Oswald); Tucson, Arizona, USA (James Clay). \par

Helmut Karzel was editor for several mathematical journals: \emph{Abhandlungen
aus dem Mathematischen Seminar der Universit\"{a}t Hamburg}, \emph{Jahresberichte
der Deutschen Mathematiker Vereinigung}, \emph{Journal of Geometry} (founding
editor and, from 2017, honorary editor), \emph{Mitteilungen der Mathematischen
Gesellschaft Hamburg}, \emph{The Nepali Mathematical Sciences Report} and
\emph{Results in Mathematics}. Since 1959, Karzel was member of the
\emph{Mathematische Gesellschaft in Hamburg}. He served as ``Jahrverwalter''
(chairman) of this society from 1966 to 1971 and was awarded an honorary
membership in 1978. In appreciation of his outstanding services to mathematics,
Karzel received an honorary doctorate from the University of Hamburg in the
year 1993. \par

Karzel inspired highly gifted young people for mathematics and encouraged them
to do own research work. In Hamburg, Karlsruhe, Hannover and M\"{u}nchen he guided
a large number of students through their theses. As a result, 30 doctorates
arose under his supervision. Several of his students later became University
Professors. \par

A list of Helmut Karzel's scientific publications, which comprises more than
170 entries, can be found at the end of this article. His scientific work is
distinguished by impressive ideas and proves his widespread interests, with
algebra, geometry and foundations of geometry playing a central role. In the
following, we will only briefly discuss a small selection of his contributions.
\par

Since the beginning of the 20th century, theorems of elementary geometry and
absolute geometry have increasingly been proven by the use of reflections. In
1943, Arnold Schmidt found a reflection-theoretical axiom system that allows
for the foundation of all classical absolute and elliptic planes. The three
reflections theorem plays an essential role in this approach. Friedrich
Bachmann gave a reduced version of Schmidt's axiom system in 1951. A fine
detailed justification and a large collection of models can be found in
Bachmann's book \emph{Aufbau der Geometrie aus dem Spiegelungsbegriff}, first
published in 1959. \par

A decisive new impetus was given to reflection geometry in 1954, when Emanuel
Sperner examined under which minimal conditions Desargues' theorem can be
proven in a group theoretic setting. Sperner based his approach on a group $G$
with a system, say $J$, of involutory generators, which are to be viewed as
``lines''. He assumed essentially only the validity of the following
\emph{general three reflections theorem}: If for five elements $A, B, X, Y, Z
\in J$, with $A\neq B$, each of the products $ABX, ABY, ABZ$ is an involution,
then $XYZ\in J$. From 1954 onwards, Karzel developed further Sperner's axiom
system in a series of papers [2, 3, 6, 7, 8, 9]. The corresponding geometries
fall into two classes:
\begin{enumerate}
\item The (previously known) \emph{regular absolute planes}. Here different
    perpendiculars to a fixed line have different pedals.

\item The so-called \emph{Lotkern-Geometrien}. These are planes in which
    all perpendiculars to a fixed line $G$ share the same pedal. Together
    with the line $G$, all perpendiculars to $G$ form a pencil, the
    \emph{Lotkern} (perpendicular pencil) of $G$.
\end{enumerate}
The geometries of the second class turned out to be new. Karzel was the first
to recognise their quite unusual properties and investigated the subject more
closely. Among the \emph{Lotkern-Geometrien} are the so-called
\emph{Zentrums-Geometrien}, where the system $J$ of generators contains an
element from the centre of the underlying group $G$. From an algebraic point of
view, \emph{Lotkern-Geometrien} arise from subgroups of orthogonal groups over
fields of characteristic two. A common foundation of the regular and the
\emph{Lotkern-Geometrien} with the help of so-called kinematic spaces (see
below) was given by Karzel in his 1963 lecture \emph{Gruppentheoretische
Begr\"{u}ndung metrischer Geometrien} (Group theoretical foundation of metric
geometries); cf.\ his book [B5], co-authored with G\"{u}nter Graumann. \par

The classical example of a kinematic space is due to Wilhelm Blaschke and,
independently, Josef Gr\"{u}nwald. In two articles, published 1911, they
established a bijection between the group of proper motions (rotations and
translations) of the Euclidean plane and a particular set of the points in the
three-dimensional real projective space. In 1934, Erich Podehl and Kurt
Reidemeister founded the elliptic plane with the help of the associated
kinematic space, which in this case comprises all points of a projective space.
Reinhold Baer extended this principle by assigning a geometry to a group $G$ in
such a way that its elements represent both points and hyperplanes. Given
$\alpha,\beta\in G$, a ``point'' $(\alpha)$ is incident with a ``hyperplane''
$\langle\beta\rangle$ if $\alpha \beta $ is involutory. The space defined in
this way is a projective space if and only if $G$ is an elliptic motion group.
\par

Karzel picked up Baer's ideas and generalised them, in collaboration with Erich
Ellers, by considering a group $G$ (with unit element $1$) together with a
distinguished subset $D$ of $G$ that is invariant under all inner automorphisms
of $G$ and such that $\xi^2=1$ for all $\xi\in D$. Then a geometry $D(G)$ is
defined as follows: For $\alpha, \beta \in G$ the incidence of $(\alpha)$ and
$\langle \beta\rangle$ means $\alpha \beta \in D$. Karzel and Ellers completely
classified the corresponding geometries: $D(G)$ is either a projective space of
dimension $3$ and $G$ is isomorphic to the motion group of an elliptic plane,
or $D(G)$ is a so-called \emph{involutory geometry} of dimension $2n-1$, which
can be represented in terms of a Clifford algebra over a field with
characteristic $2$ (cf.\ [11, 12]). \par

Let us take a closer look at the geometry $D(G)$ from above. For all $a\in G$,
the left translations $a_\ell\colon x \to ax$ and the right translations $ a_r
\colon x \to xa$ are automorphisms of $D(G)$ that act regularly on the set of
points. Starting from this observation, Karzel developed the concept of an
\emph{incidence group}, a group with a ``compatible'' geometric structure,
comparable to the concept of a topological group. To be more precise, a group
$G$ is called an \emph{incidence group} if the elements of $G$ are the
``points'' of a geometry (with ``lines'' understood as subsets of $G$) and the
left translations in $G$ are automorphisms of the geometry; see the surveys
[24] and [31], the latter being joint work with Irene Pieper. Motivated by
previous examples, Karzel also coined the notion of an abstract \emph{kinematic
space} [40] as an incidence group satisfying two additional conditions: (i) all
right translations in $G$ are automorphisms of the geometry; (ii) lines through
the neutral element of $G$ are subgroups of $G$. \par

In numerous publications and several dissertations under his guidance, Karzel
advanced the algebraic description of incidence groups in terms of
representation theorems. \emph{Normal nearfields}, i.e. nearfields which are
also left vector spaces over a normal sub-nearfield, play a crucial role in
some of these results. \par

Kinematic spaces permit two \emph{parallelisms} by considering the left (or
right) cosets the lines passing the neutral element (of the underlying group)
as classes of parallel lines. It turned out that kinematic spaces can be
characterised as incidence geometries with two parallelisms; see [41, 44, 45,
86, 88, 89, 93] (with Erich Ellers, Hans-Joachim Kroll, Carl J.~Maxson, Kay
S\"{o}rensen). \par

In the above mentioned nearfields, the validity of one distributive law is
removed as a generalisation to fields. Leonard E.~Dickson gave examples of
nearfields that are not fields. Hans Zassenhaus extended Dickson's construction
method and showed that in this way, up to seven exceptions, all finite
nearfields can be obtained. Karzel axiomatised this construction method
(Dickson's process) in [21] and characterised those nearfields that can be
obtained in this way, known as \emph{Dickson's nearfields}. \par

Given a nearfield $F$ the group
\begin{equation*}
    T_2(F) := \{ \tau_{a,b}\colon F \to F;\, x \to a+ bx \mid a \in F, b \in F^*\}
\end{equation*}
of all affine translations operates sharply $2$-transitively on $F$.
Conversely, if $T$ is a sharply $2$-transitive permutation group on a finite
set $F$, then, according to Robert D.~Carmichael, two operations can be defined
on $F$ such that $F$ is made into a nearfield and $T$ is isomorphic to the
group $T_2(F)$. In the infinite case, the situation is more intricate. Karzel
introduced two operations on $F$, translated the double transitivity into
algebraic axioms and obtained a so-called \emph{neardomain}. In a neardomain
$F$ the additive structure is in general not a group, but only a loop with
additional properties, a so-called \emph{K-loop} (in particular, there is an
additive automorphism $\delta_{a,b}$ with $a+(b+x)= (a+b)+\delta_{a,b} (x)$ for
any elements $a,b \in F$). Compared to other concepts, neardomains can be
characterised by the fact that isomorphic permutation groups lead to isomorphic
neardomains and vice versa. \par

The existence of proper K-loops was an open problem for a long time.
Surprisingly, the first proper example of a K-loop was found in the context of
mathematical physics. Abraham A.~Ungar investigated in 1988 the relativistic
velocity addition $\oplus$ on $\mathbb{R}^3_c= \{ x \in \mathbb{R}^3: |x| < c
\}$, which is neither commutative nor associative. However, Ungar was able to
prove that $\oplus$ makes $\mathbb{R}^3_c$ into a K-loop. Investigations on
K-loops and relativistic velocity addition are in focus of Karzel's papers
[100, 104, 106] (with Bokhee Im, Heinrich Wefelscheid). \par

Last but not least, let us take a glimpse at the remaining areas of Karzel's
work. These include, among others, \emph{circle geometries}, intensively worked
on by Karzel in [36, 37, 38, 39, 48, 66, 68, 91] (with Werner Heise,
Hans-Joachim Kroll, Helmut M\"{a}urer, Rotraut Stanik, Heinz W\"{a}hling), \emph{coding
theory} (in particular applications of circle geometry in coding theory) [90,
91, 99] (with Alan Oswald, Carl J.~Maxson), questions of the \emph{order
relation} in algebra and geometry [1, 4, 5, 13, 14, 16, 30, 48, 94] (with
Hanfried Lenz, Rotraut Stanik, Heinz W\"{a}hling) and the foundation of
\emph{metric planes} [47, 55, 56, 57, 61, 65, 77] (with G\"{u}nter Kist, Rotraut
Stanik, Monika K\"{o}nig).

\section*{Articles}

\begin{enumerate}
\item %\label{}
    H.~Karzel, Erzeugbare Ordnungsfunktionen. \emph{Math.\ Ann.}
    \textbf{127} (1954), 228--242.

\item %\label{}
    H.~Karzel, Ein Axiomensystem der absoluten Geometrie. \emph{Arch.\
    Math.\ (Basel)} \textbf{6} (1954), 66--76.

\item %\label{}
    H.~Karzel, Verallgemeinerte absolute Geometrien und Lotkerngeometrien.
    \emph{Arch.\ Math.\ (Basel)} \textbf{6} (1955), 284--295.

\item %\label{}
    H.~Karzel, Ordnungsfunktionen in nichtdesarguesschen projektiven
    Geometrien. \emph{Math.\ Z.} \textbf{62} (1955), 268--291.

\item %\label{}
    H.~Karzel, \"{U}ber eine Anordnungsbeziehung am Dreieck. \emph{Math.\ Z.}
    \textbf{64} (1956), 131--137.

\item %\label{}
    H.~Karzel, Kennzeichnung der Gruppe der gebrochen-linearen
    Transformationen \"{u}ber einem K\"{o}rper der Charakteristik $2$. \emph{Abh.\
    Math.\ Sem.\ Univ.\ Hamburg} \textbf{22} (1958), 1--8.

\item %\label{}
    H.~Karzel, Spiegelungsgeometrien mit echtem Zentrum. \emph{Arch.\
    Math.\ (Basel)} \textbf{9} (1958), 140--146.

\item %\label{}
    H.~Karzel, Zentrumsgeometrien und elliptische Lotkerngeometrien.
    \emph{Arch.\ Math.\ (Basel)} \textbf{9} (1958), 455--464.

\item %\label{}
    H.~Karzel, Quadratische Formen von Geometrien der Charakteristik $2$.
    \emph{Abh.\ Math.\ Sem.\ Univ.\ Hamburg} \textbf{23} (1959), 144--162.

\item %\label{}
    H.~Karzel, Wandlungen des Begriffs der projektiven Geometrie.
    \emph{Mitt.\ Math.\ Ges.\ Hamburg} \textbf{9} (1959), 42--48.

\item %\label{}
    H.~Karzel, Verallgemeinerte elliptische Geometrien und ihre
    Gruppenr\"{a}ume. \emph{Abh.\ Math.\ Sem.\ Univ.\ Hamburg} \textbf{24}
    (1960), 167--188.

\item %\label{}
    E.~Ellers, H.~Karzel, Involutorische Geometrien. \emph{Abh.\ Math.\
    Sem.\ Univ.\ Hamburg} \textbf{25} (1961), 93--104.

\item %\label{}
    H.~Karzel, H.~Lenz, \"{U}ber Hilbertsche und Spernersche Anordnung.
    \emph{Abh.\ Math.\ Sem.\ Univ.\ Hamburg} \textbf{25} (1961), 82--88.

\item %\label{}
    H.~Karzel, Anordnungsfragen in tern\"{a}ren Ringen und allgemeinen
    projektiven und affinen Ebenen. In: \emph{Algebraical and Topological
    Foundations of Geometry (Proceedings of a Colloquium held in Utrecht,
    August 1959)}, 71--86, Pergamon, New York 1962.

\item %\label{}
    H.~Karzel, Kommutative Inzidenzgruppen. \emph{Arch.\ Math.\ (Basel)}
    \textbf{13} (1962), 535--538.

\item %\label{}
    H.~Karzel, Zur Fortsetzung affiner Ordnungsfunktionen. \emph{Abh.\
    Math.\ Sem.\ Univ.\ Hamburg} \textbf{26} (1963), 17--22.

\item %\label{}
    E.~Ellers, H.~Karzel, Kennzeichnung elliptischer Gruppenr\"{a}ume.
    \emph{Abh.\ Math.\ Sem.\ Univ.\ Hamburg} \textbf{26} (1963), 55--77.

\item %\label{}
    H.~Karzel, Ebene Inzidenzgruppen. \emph{Arch.\ Math.\ (Basel)}
    \textbf{15} (1964), 10--17.

\item %\label{}
    E.~Ellers, H.~Karzel, Endliche Inzidenzgruppen. \emph{Abh.\ Math.\
    Sem.\ Univ.\ Hamburg} \textbf{27} (1964), 250--264.

\item %\label{}
    H.~Karzel, Beziehungen zwischen topologischen Inzidenzgruppen und
    topologischen Fastk\"{o}rpern. In: \emph{Simposio di Topologia
    (Celebrazioni archimedee del secolo XX, Messina, 1964)}, 75--84,
    Edizioni Oderisi, Gubbio 1965.

\item %\label{}
    H.~Karzel, Unendliche Dicksonsche Fastk\"{o}rper. \emph{Arch.\ Math.\
    (Basel)} \textbf{16} (1965), 247--256.

\item %\label{}
    H.~Karzel, Projektive R\"{a}ume mit einer kommutativen transitiven
    Kollineationsgruppe. \emph{Math.\ Z.} \textbf{87} (1965), 74--77.

\item %\label{}
    H.~Karzel, Normale Fastk\"{o}rper mit kommutativer Inzidenzgruppe.
    \emph{Abh.\ Math.\ Sem.\ Univ.\ Hamburg} \textbf{28} (1965), 124--132.

\item %\label{}
    H.~Karzel, Bericht \"{u}ber projektive Inzidenzgruppen. \emph{Jber.\
    Deutsch.\ Math.-Verein.} \textbf{67} (1965), 58--92.

\item %\label{}
    H.~Karzel, Zweiseitige Inzidenzgruppen. \emph{Abh.\ Math.\ Sem.\ Univ.\
    Hamburg} \textbf{29} (1965), 118--136.

\item\label{pub:ellers+karzel1967} %
    E.~Ellers, H.~Karzel, Die klassische euklidische und hyperbolische
    Geometrie. In: H.~Behnke, F.~Bachmann, K.~Fladt (eds.), \emph{Grundz\"{u}ge
    der Mathematik, Bd.~II, Geometrie, Teil~A: Grundlagen der Geometrie,
    Elementargeometrie}, 187--213, Vandenhoeck \& Ruprecht, G\"{o}ttigen 1967.

\item %\label{}
    H.~Karzel, H.~Mei{\ss}ner, Geschlitzte Inzidenzgruppen und normale
    Fastmoduln. \emph{Abh.\ Math.\ Sem.\ Univ.\ Hamburg} \textbf{31}
    (1967), 69--88.

\item %\label{}
    H.~Karzel, Projective planes with a commutative and transitive
    collineation group. In: R.~Sandler (ed.), \emph{Proceedings of the
    Projective Geometry Conference at the University of Illinois, Summer
    1967}, 80--82, Dept.\ of Mathematics, University of Illinois, Chicago
    1967.

\item %\label{}
    H.~Karzel, Zusammenh\"{a}nge zwischen Fastbereichen, scharf zweifach
    transitiven Permutationsgruppen und $2$-Strukturen mit Rechtecksaxiom.
    \emph{Abh.\ Math.\ Sem.\ Univ.\ Hamburg} \textbf{32} (1968), 191--206.

\item %\label{}
    H.~Karzel, Konvexit\"{a}t in halbgeordneten projektiven und affinen R\"{a}umen.
    \emph{Abh.\ Math.\ Sem.\ Univ.\ Hamburg} \textbf{33} (1969), 231--242.

\item %\label{}
    H.~Karzel, I.~Pieper, Bericht \"{u}ber geschlitzte Inzidenzgruppen.
    \emph{Jber.\ Deutsch.\ Math.-Verein.} \textbf{72} (1970), 70--114.

\item %\label{}
    H.~Karzel, Spiegelungsgruppen und absolute Gruppenr\"{a}ume. \emph{Abh.\
    Math.\ Sem.\ Univ.\ Hamburg} \textbf{35} (1971), 141--163.

\item %\label{}
    H.~Karzel, K.~S\"{o}rensen, Die lokalen S\"{a}tze von Pappus und Pascal.
    \emph{Mitt.\ Math.\ Ges.\ Hamburg} \textbf{10} (1971), 28--55.

\item %\label{}
    H.~Karzel, K.~S\"{o}rensen, Projektive Ebenen mit einem pascalschen Oval.
    \emph{Abh.\ Math.\ Sem.\ Univ.\ Hamburg} \textbf{36} (1971), 123--125.

\item %\label{}
    E.~Ellers, H.~Karzel, Involutory incidence spaces. \emph{J.\ Geom.}
    \textbf{1} (1971), 117--126.

\item %\label{}
    W.~Heise, H.~Karzel, Eine Charakterisierung der ovoidalen
    Kettengeometrien. \emph{J.\ Geom.} \textbf{2} (1972), 69--74.

\item %\label{}
    W.~Heise, H.~Karzel, Laguerre- und Minkowski-{$m$}-Strukturen.
    \emph{Rend.\ Istit.\ Mat.\ Univ.\ Trieste} \textbf{4} (1972), 139--147.

\item %\label{}
    W.~Heise, H.~Karzel, Vollkommen fanosche Minkowski-Ebenen. \emph{J.\
    Geom.} \textbf{3} (1973), 21--29.

\item %\label{}
    W.~Heise, H.~Karzel, Symmetrische Minkowski-Ebenen. \emph{J.\ Geom.}
    \textbf{3} (1973), 5--20.

\item %\label{}
    H.~Karzel, Kinematic spaces. In: \emph{Symposia Mathematica, Vol.~XI
    (Convegno di Geometria, Istituto Nazionale di Alta Matematica, Roma,
    22--25 maggio 1972)}, 413--439, Academic Press, London 1973.

\item %\label{}
    H.~Karzel, H.-J.~Kroll, K.~S\"{o}rensen, Invariante Gruppenpartitionen und
    Doppelr\"{a}ume. \emph{J.\ Reine Angew.\ Math.} \textbf{262-263} (1973),
    153--157.

\item %\label{}
    H.~Karzel, Endliche $2$-gelochte Ebenen. \emph{Abh.\ Math.\ Sem.\
    Univ.\ Hamburg} \textbf{40} (1974), 187--196.

\item %\label{}
    H.~Karzel, Kinematische Algebren und ihre geometrischen Ableitungen.
    \emph{Abh.\ Math.\ Sem.\ Univ.\ Hamburg} \textbf{41} (1974), 158--171.

\item %\label{}
    H.~Karzel, H.-J.~Kroll, K.~S\"{o}rensen, Projektive Doppelr\"{a}ume.
    \emph{Arch.\ Math.\ (Basel)} \textbf{25} (1974), 206--209.

\item %\label{}
    E.~Ellers, H.~Karzel, The classical Euclidean and the classical
    hyperbolic geometry. In: H.~Behnke, F.~Bachmann, K.~Fladt, H.~Kunle
    (eds.), \emph{Fundamentals of mathematics, Vol. II: Geometry
    (Translated from the second German edition by S.~H.~Gould)}, 174--197,
    MIT Press, Cambridge, Mass.\ 1974. [Cf.~\ref{pub:ellers+karzel1967}.]

\item %\label{}
    H.~Karzel, H.-J.~Kroll, Eine inzidenzgeometrische Kennzeichnung
    projektiver kinematischer R\"{a}ume. \emph{Arch.\ Math.\ (Basel)}
    \textbf{26} (1975), 107--112.

\item %\label{}
    H.~Karzel, Fanosche metrische affine Ebenen. \emph{Abh.\ Math.\ Sem.\
    Univ.\ Hamburg} \textbf{43} (1975), 166--178.

\item %\label{}
    H.~Karzel, R.~Stanik, H.~W\"{a}hling, Zum Anordnungsbegriff in affinen
    Geometrien. \emph{Abh.\ Math.\ Sem.\ Univ.\ Hamburg} \textbf{44}
    (1976), 24--31.

\item %\label{}
    H.~Karzel, Some recent results on incidence groups. In: P.~Scherk
    (ed.), \emph{Foundations of Geometry: Selected Proceedings of a
    Conference (University of Toronto, July 17 to August 18, 1974)},
    114--143, University of Toronto Press, Toronto 1976.

\item %
    H.~Karzel, Kongruenzen in projektiven R\"{a}umen und ihre $2$-dimensionalen
    Ableitungen. In: \emph{Kolloquium \"{u}ber Geometrie 1975, II}, 1--10,
    Institut f\"{u}r Mathematik TU Hannover, Bericht Nr.\ \textbf{45}, Hannover
    1976.

\item %\label{}
    H.~Karzel, H.-J.~Kroll, Gruppen von Projektivit\"{a}ten in Zwei- und
    Hyperbelstrukturen. In: \emph{Lenz Festband (Hanfried Lenz zu seinem
    60.~Geburtstage)}, 125--134, Fachbereich Mathematik FU Berlin, Preprint
    Nr.\ \textbf{9}, Berlin 1976.

\item %\label{}
    H.~Karzel, Symmetrische Permutationsmengen. \emph{Aequationes Math.}
    \textbf{17} (1978), 83--90.

\item %\label{}
    H.~Karzel, G.~Kist, H.-J.~Kroll, Burau-Geometrien. \emph{Results Math.}
    \textbf{2} (1979), 88--104.

\item %\label{}
    H.~Karzel, H.-J.~Kroll, Zur projektiven Einbettung von Inzidenzr\"{a}umen
    mit Eigentlichkeitsbereich. \emph{Abh.\ Math.\ Sem.\ Univ.\ Hamburg}
    \textbf{49} (1979), 82--94.

\item %\label{}
    H.~Karzel, G.~Kist, Zur Begr\"{u}ndung metrischer affiner Ebenen.
    \emph{Abh.\ Math.\ Sem.\ Univ.\ Hamburg} \textbf{49} (1979), 234--236.

\item %\label{}
    H.~Karzel, R.~Stanik, Metrische affine Ebenen. \emph{Abh.\ Math.\ Sem.\
    Univ.\ Hamburg} \textbf{49} (1979), 237--243.

\item %\label{}
    H.~Karzel, R.~Stanik, Rechtseitebenen und ihre Darstellung durch
    Integrit\"{a}tssysteme. \emph{Mitt.\ Math.\ Ges.\ Hamburg} \textbf{10}
    (1979), 531--551.

\item %\label{}
    H.~Karzel, Zyklisch geordnete Gruppen. \emph{Mitt.\ Math.\ Ges.\
    Hamburg} \textbf{10} (1979), 523--529.

\item H.~Karzel, Einige neuere Beitr\"{a}ge zur Theorie der Inzidenzgruppen.
    In: Geometrie Seminar, 9.\ bis 13.~Mai 1977, Aristoteles Universit\"{a}t,
    Thessaloniki, \emph{J.\ Geom} \textbf{13} (1979), 12--15.

\item %\label{}
    H.~Karzel, G.~Kist, Some applications of near-fields. \emph{Proc.\
    Edinburgh Math.\ Soc.\ (2)} \textbf{23} (1980), 129--139.

\item %\label{}
    H.~Karzel, Rectangular spaces. In: R.~Artzy, I.~Vaisman (eds.),
    \emph{Geometry and Differential Geometry (Proceedings of a Conference
    Held at the University of Haifa, Israel, March 18--23, 1979)}, volume
    792 of \emph{Lecture Notes in Math.}, 79--91, Springer, Berlin 1980.

\item %\label{}
    H.~Karzel, H.-J.~Kroll, Zur Inzidenzstruktur kinematischer R\"{a}ume und
    absoluter Ebenen. \emph{Beitr\"{a}ge zur Geometrie und Algebra Nr.\ 6}, TU
    M\"{u}nchen, Math.\ Inst., \textbf{M8010} (1980), 42--61.

\item %\label{}
    H.~Karzel, Emanuel Sperner als Begr\"{u}nder einer neuen Anordnungstheorie.
    In: H.~Zeitler (ed.), \emph{Gedenkkolloquium f\"{u}r Dr.~Dr.~h.~c. Emanuel
    Sperner}, 21--37, Universit\"{a}t Bayreuth, Bayreuth 1980.

\item %\label{}
    H.~Karzel, Spazi cinematici e geometria di riflessione. \emph{Quad.\
    Sem.\ Geometrie Combinatorie} n.~\textbf{30} (Settembre 1980),
    Dipartimento di Matematica, Universit\`{a} di Roma ``La Sapienza'', 1--115.

\item %\label{}
    H.~Karzel, M.~K\"{o}nig, Affine Einbettung absoluter R\"{a}ume beliebiger
    Dimension. In: P.~L.~Butzer, F.~Feh\'{e}r (eds.), \emph{E.~B.~Christoffel:
    The Influence of His Work on Mathematics and the Physical Sciences},
    657--670, Birkh\"{a}user, Basel 1981.

\item %\label{}
    H.~Karzel, H.-J.~Kroll, Perspectivities in circle geometries. In:
    P.~Plaumann, K.~Strambach (eds.), \emph{Geometry~-- von Staudt's point
    of view (Proceedings of the NATO Advanced Study Institute held at Bad
    Windsheim, West Germany, July~21--August~1, 1980)}, volume~70 of
    \emph{NATO Adv.\ Study Inst.\ Ser.~C: Math.\ Phys.\ Sci.}, 51--99,
    Reidel, Dordrecht 1981.

\item %\label{}
    H.~Karzel, The projectivity groups of quadratic sets and their
    representations by neardomains and nearfields. In: G.~Ferrero,
    C.~Ferrero Cotti (eds.), \emph{Proceedings of Conference on Near-Rings
    and Near-Fields (San Benedetto del Tronto, September 13--19, 1981)},
    95-100, Universit\`{a} degli Studi di Parma, Parma, 1982.

\item %\label{}
    H.~Karzel, H.~M\"{a}urer, Eine Kennzeichnung miquelscher M\"{o}biusebenen durch
    eine Eigenschaft der Kreisspiegelungen. \emph{Results Math.} \textbf{5}
    (1982), 52--56.

\item %\label{}
    H.~Karzel, Erinnerungen an Emanuel Sperner aus den Jahren 1948--1968
    und Emanuel Sperners Beitr\"{a}ge zur metrischen Geometrie und ihre
    Bedeutung f\"{u}r die Entwicklung der Geometrie. \emph{Mitt.\ Math.\ Ges.\
    Hamburg} \textbf{11} (1983), 217--231.

\item %\label{}
    H.~Karzel, M.~Marchi, The projectivity groups of ovals and of quadratic
    sets. In: A.~Barlotti, P.~V.~Ceccherini, G.~Tallini (eds.),
    \emph{Combinatorics~'81, In Honour of Beniamino Segre (Proceedings of
    the International Conference on Combinatorial Geometries and their
    Applications, Rome, June 7--12, 1981, Annals of Discrete Math.~18)},
    volume~78 of \emph{North-Holland Math.\ Stud.}, 519--533,
    North-Holland, Amsterdam 1983.

\item %\label{}
    H.~Karzel, M.~Marchi, Plane fibered incidence groups. \emph{J.\ Geom.}
    \textbf{20} (1983), 192--201.

\item %\label{}
    H.~Karzel, Affine incidence groups. In: A.~Barlotti, M.~Marchi,
    G.~Tallini (eds.), \emph{Atti del Convegno ``Geometria Combinatoria e
    di Incidenza: Fondamenti e Applicazioni'' (La Mendola, 4--11 Luglio
    1982)}, volume~7 of \emph{Rend.\ Sem.\ Mat.\ Brescia}, 409--425, Vita e
    Pensiero, Milano 1984.

\item %\label{}
    H.~Karzel, C.~J.~Maxson, Fibered groups with non-trivial centers.
    \emph{Results Math.} \textbf{7} (1984), 192--208.

\item %\label{}
    H.~Karzel, C.~J.~Maxson, Kinematic spaces with dilatations. \emph{J.\
    Geom.} \textbf{22} (1984), 196--201.

\item %\label{}
    H.~Karzel, Fastvektorr\"{a}ume, unvollst\"{a}ndige Fastk\"{o}rper und ihre
    abgeleiteten geometrischen Strukturen. \emph{Mitt.\ Math.\ Sem.\
    Giessen} \textbf{166} (1984), 127--139.

\item %\label{}
    H.~Karzel, G.~Kist, Determination of all near vector spaces with
    projective and affine fibrations. \emph{J.\ Geom.} \textbf{23} (1984),
    124--127.

\item %\label{}
    H.~Karzel, Zur Begr\"{u}ndung euklidischer R\"{a}ume. \emph{Mitt.\ Math.\ Ges.\
    Hamburg} \textbf{11} (1985), 355--368.

\item %\label{}
    H.~Karzel, G.~Kist, Kinematic algebras and their geometries. In:
    R.~Kaya, P.~Plaumann, K.~Strambach (eds.), \emph{Rings and Geometry
    (Proceedings of the NATO Advanced Study Institute, Istanbul, Turkey,
    September 2--14, 1984)}, volume 160 of \emph{NATO Adv.\ Sci.\ Inst.\
    Ser.~C Math.\ Phys.\ Sci.}, 437--509, Reidel, Dordrecht 1985.

\item %\label{}
    H.~Karzel, G.~Kist, Finite porous incidence groups and s-homogeneous
    $2$-sets. In: \emph{Symposia Mathematica, Vol.~XXVIII (Combinatorica,
    Convegno, Istituto Nazionale di Alta Matematica, Roma, 23--26 maggio
    1983)}, 89--111, Academic Press, London 1986.

\item %\label{}
    H.~Karzel, C.J.~Maxson, Fibered $p$-groups. \emph{Abh.\ Math.\ Sem.\
    Univ.\ Hamburg} \textbf{56} (1986), 1--9.

\item %\label{}
    H.~Karzel, C.~J.~Maxson, G.~Pilz, Kernels of covered groups.
    \emph{Results Math.} \textbf{9} (1986), 70--81.

\item %\label{}
    J.~R.~Clay, H.~Karzel, Tactical configurations derived from groups
    having a group of fixed point free automorphisms. \emph{J.\ Geom.}
    \textbf{27} (1986), 60--68.

\item %\label{}
    H.~Karzel, \"{U}ber einen Fundamentalsatz der synthetischen algebraischen
    Geometrie von W.~Burau und H.~Timmermann. \emph{J.\ Geom.} \textbf{28}
    (1987), 86--101.

\item %\label{}
    H.~Karzel, Couplings and derived structures. In: G.~Betsch (ed.),
    \emph{Near-rings and Near-fields (Proceedings of a Conference held at
    the University of T\"{u}bingen, F.R.G., 4--10 August, 1985)}, volume 137 of
    \emph{North-Holland Math.\ Stud.}, 133--143, North-Holland, Amsterdam
    1987.

\item %\label{}
    J.~R.~Clay, H.~Karzel, On the equality of blocks of a group of fixed
    point free automorphisms. \emph{Results Math.} \textbf{13} (1988),
    33--46.

\item %\label{}
    H.~Karzel, Double spaces of order $3$ and D-loops of exponent $3$. In:
    R.~Mlitz (ed.), \emph{General Algebra 1988 (Proceedings of the
    International Conference, Held in Memory of Wilfried N\"{o}bauer, Krems,
    Austria, August 21--27, 1988)}, 115--127, North-Holland, Amsterdam
    1990.

\item %\label{}
    H.~Karzel, M.~Marchi, Axiomatic characterization of weak projective
    spaces. \emph{Riv.\ Mat.\ Pura Appl.} \textbf{3} (1988), 99--105
    (1989).

\item %\label{}
    H.~Karzel, Porous double spaces. \emph{J.\ Geom.} \textbf{34} (1989),
    80--104.

\item %\label{}
    H.~Karzel, Affine double spaces of order $3$. \emph{Results Math.}
    \textbf{15} (1989), 75--80.

\item %\label{}
    H.~Karzel, A.~Oswald, Near-rings (MDS)- and Laguerre codes. \emph{J.\
    Geom.} \textbf{37} (1990), 105--117.

\item %\label{}
    H.~Karzel, Circle geometry and its application to code theory. In:
    G.~Longo, M.~Marchi, A.~Sgarro (eds.), \emph{Geometries, Codes and
    Cryptography (Papers from the school held in Udine, 1989)}, volume 313
    of \emph{CISM Courses and Lect.}, 25--75, Springer-Verlag, Wien 1990.

\item %\label{}
    H.~Karzel, B.~Reinmiedl, Pseudo-subaffine double spaces. In:
    A.~Barlotti, G.~Lunardon, F.~Mazzocca, N.~Melone, D.~Olanda, A.~Pasini,
    G.~Tallini (eds.), \emph{Combinatorics '88, Vol.~2 (Proceedings of the
    International Conference on Incidence Geometries and Combinatorial
    Structures held in Ravello, May 23--28, 1988)}, 129--139, Mediterranean
    Press, Rende 1991.

\item %\label{}
    H.~Karzel, C.~J.~Maxson, Archimedeisation of some ordered geometric
    structures which are related to kinematic spaces. \emph{Results Math.}
    \textbf{19} (1991), 290--318.

\item %\label{}
    H.~Karzel, Finite reflection groups and their corresponding structures.
    In: A.~Barlotti, A.~Bichara, P.~V.~Ceccherini, G.~Tallini (eds.),
    \emph{Combinatorics '90. Recent trends and applications (Proceedings of
    the International Conference held in Gaeta, May 20--27, 1990)},
    volume~52 of \emph{Ann.\ Discrete Math.}, 317--336, North-Holland,
    Amsterdam 1992.

\item %\label{}
    H.~Karzel, M.~J.~Thomsen, Near-rings, generalizations, near-rings with
    regular elements and applications, a report. In: G.~Pilz (ed.),
    \emph{Contributions to General Algebra 8 (Proceedings of the Conference
    on Near-rings and Near-fields held in Linz, July 14--20, 1991)},
    91--110, H\"{o}lder-Pichler-Tempsky, Wien 1992.

\item %\label{}
    H.~Karzel, M.~Marchi, S.~Pianta, On commutativity in point-reflection
    geometries. \emph{J.\ Geom.} \textbf{44} (1992), 102--106.

\item %\label{}
    H.~Karzel, Geometrie affini, geometrie metriche, e loro trasformazioni.
    In: \emph{Nuova Secondaria 1992 n.~10}, 44--48, Editrice La Scuola,
    Brescia 1992.

\item %\label{}
    H.~Karzel, Kinematic structures of generalized hyperbolic spaces. In:
    A.~Barlotti, M.~Gionfriddo, G.~Tallini (eds.), \emph{Combinatorics '92
    (Lectures from the Third Combinatorial Conference held in Catania,
    September 6--13, 1992.)}, \emph{Matematiche (Catania)} \textbf{47}
    (1992), 259--279.

\item %\label{}
    H.~Karzel, C.~J.~Maxson, Affine {MDS}-codes on groups. \emph{J.\ Geom.}
    \textbf{47} (1993), 65--76.

\item %\label{}
    H.~Karzel, H.~Wefelscheid, Groups with an involutory antiautomorphism
    and K-loops; application to space-time-world and hyperbolic geometry~I.
    \emph{Results Math.} \textbf{23} (1993), 338--354.

\item %\label{}
    H.~Karzel, S.~Pianta, R.~Stanik, Generalized Euclidean and elliptic
    geometries, their connections and automorphism groups. \emph{J.\ Geom.}
    \textbf{48} (1993), 109--143.

\item %\label{}
    H.~Karzel, The Lorentz group and the hyperbolic geometry.
    \emph{Beitr\"{a}ge zur Geometrie und Algebra Nr.\ 24}, TU M\"{u}nchen, Math.\
    Inst., \textbf{M9315} (1993), 5--9.

\item %\label{}
    H.~Karzel, A.~Oswald, Lecture notes on algebras of $2\times 2$-matrices
    over quadratic field extension and their geometric derivations.~I:
    Fundamental properties. \emph{Beitr\"{a}ge zur Geometrie und Algebra Nr.\
    25}, TU M\"{u}nchen, Math.\ Inst., \textbf{M9317} (1993), 36~pp.

\item %\label{}
    B.~Im, H.~Karzel, Determination of the automorphism group of a
    hyperbolic K-loop. \emph{J.\ Geom.} \textbf{49} (1994), 96--105.

\item %\label{}
    H.~Karzel, A.~Konrad, Eigenschaften angeordneter R\"{a}ume mit
    hyperbolischer Inzidenzstruktur~I. \emph{Beitr\"{a}ge zur Geometrie und
    Algebra Nr.\ 28}, TU M\"{u}nchen, Math.\ Inst., \textbf{M9415} (1994),
    27--36.

\item %\label{}
    H.~Karzel, Raum-Zeit-Welt und hyperbolische Geometrie (Ausarbeitung der
    von Prof.~Dr.~Dr.~h.~c.\ Helmut Karzel im SS~1992 und WS~1992/93 an der
    Technischen Universit\"{a}t M\"{u}nchen gehaltenen Vorlesungen, ausgearbeitet
    von Angelika Konrad). \emph{Beitr\"{a}ge zur Geometrie und Algebra Nr.\
    29}, TU M\"{u}nchen, Math.\ Inst., \textbf{M9412} (1994), 175~pp.

\item %\label{}
    H.~Karzel, A.~Konrad, A.~Kreuzer, Zur projektiven Einbettung
    angeordneter R\"{a}ume mit hyperbolischer Inzidenzstruktur. \emph{Beitr\"{a}ge
    zur Geometrie und Algebra Nr.\ 30}, TU M\"{u}nchen, Math.\ Inst.,
    \textbf{M9502} (1995), 17--27.

\item %\label{}
    H.~Karzel, A.~Konrad, A.~Kreuzer, Eigenschaften angeordneter R\"{a}ume mit
    hyperbolischer Inzidenzstruktur~II. \emph{Beitr\"{a}ge zur Geometrie und
    Algebra Nr.\ 33}, TU M\"{u}nchen, Math.\ Inst., \textbf{M9509} (1995),
    7--14.

\item H.~Karzel, Hilberts Einflu{\ss} auf die Entwicklung der Geometrie. In:
    \emph{Acta Borussica, Beitr\"{a}ge zur ost- und westpreu{\ss}ischen
    Landeskunde, Band V}, Altpreu{\ss}ische Gesellschaft f\"{u}r Wissenschaft,
    Kunst und Literatur, M\"{u}nchen 1995.

\item %\label{}
    H.~Karzel, Development of non-Euclidean geometries since Gau{\ss}. In:
    M.~Behara, R.~Fritsch, R.~G.~Lintz (eds.), \emph{Symposia Gaussiana.
    Conference~A: Mathematics and Theoretical Physics (Proceedings of the
    2nd Gauss Symposium, Munich, August 2--7, 1993)}, 397--417, de Gruyter,
    Berlin 1995.

\item %\label{}
    B.~Im, H.~Karzel, K-loops over dual numbers. \emph{Results Math.}
    \textbf{28} (1995), 75--85.

\item %\label{}
    H.~Karzel, A.~Konrad, Reflection groups and K-loops. \emph{J.\ Geom.}
    \textbf{52} (1995), 120--129.

\item %\label{}
    H.~Karzel, H.~Wefelscheid, A geometric construction of the K-loop of a
    hyperbolic space. \emph{Geom. Dedicata} \textbf{58} (1995), 227--236.

\item %\label{}
    J.~C.~Fisher, H.~Karzel, H.~Kiechle, Bundles of conics derived from
    planar projective incidence groups. \emph{J.\ Geom.} \textbf{59}
    (1997), 34--45.

\item %\label{}
    H.~Karzel, From nearrings and nearfields to K-loops. In: G.~Saad,
    M.~J.~Thomsen (eds.), \emph{Nearrings, Nearfields and K-Loops
    (Proceedings of the Conference on Nearrings and Nearfields, Hamburg,
    Germany, July 30--August 6, 1995)}, volume 426 of \emph{Math.\ Appl.},
    1--20, Kluwer Acad.\ Publ., Dordrecht 1997.

\item %\label{}
    H.~Karzel, Survey on the papers and personal memories of Giuseppe
    Tallini. \emph{Results Math.} \textbf{32} (1997), 260--271.

\item %\label{}
    E.~Gabrieli, H.~Karzel, Point-reflection geometries, geometric K-loops
    and unitary geometries. \emph{Results Math.} \textbf{32} (1997),
    66--72.

\item %\label{}
    E.~Gabrieli, H.~Karzel, Reflection geometries over loops. \emph{Results
    Math.} \textbf{32} (1997), 61--65.

\item %\label{}
    E.~Gabrieli, H.~Karzel, The reflection structures of generalized
    co-Minkowski spaces leading to K-loops. \emph{Results Math.}
    \textbf{32} (1997), 73--79.

\item %\label{}
    H.~Karzel, Geometric reflection structures and their corresponding
    loops. \emph{Rend.\ Circ.\ Mat.\ Palermo (2) Suppl.} \textbf{53}
    (1998), 119--129.

\item %\label{}
    H.~Karzel, M.~Marchi, L'orientamento nel piano: una difficile
    razionalizzazione. In: B.~D'Amore, C.~Pellegrino (eds.), \emph{Convegno
    per i sessantacinque anni di Francesco Speranza (Bologna, Dipartimento
    di matematica, sabato 11 ottobre 1997)}, 65--73, Pitagora Editrice,
    Bologna 1998.

\item %\label{}
    H.~Karzel, M.~Marchi, Regular incidence permutation sets and incidence
    quasigroups. \emph{J.\ Geom.} \textbf{63} (1998), 109--123.

\item %\label{}
    E.~Gabrieli, B.~Im, H.~Karzel, Webs related to K-loops and reflection
    structures. \emph{Abh.\ Math.\ Sem.\ Univ.\ Hamburg} \textbf{69}
    (1999), 89--102.

\item %\label{}
    B.~Im, H.~Karzel, Centralizers of certain matrices relative to the
    operation related to a K-loop. \emph{Results Math.} \textbf{36} (1999),
    69--74.

\item %\label{}
    H.~Karzel, Recent developments on absolute geometries and
    algebraization by K-loops. \emph{Discrete Math.} \textbf{208/209}
    (1999), 387--409.

\item %\label{}
    H.~Karzel, E.~Zizioli, Extension of a class of fibered loops to
    kinematic spaces. \emph{J.\ Geom.} \textbf{65} (1999), 117--129.

\item %\label{}
    H.~Karzel, H.~Wefelscheid, Werner Burau (1906--1994). \emph{Mitt.\
    Math.\ Ges.\ Hamb.} \textbf{19*} (2000), 167--183.

\item %\label{}
    B.~Alinovi, H.~Karzel, C.~Tonesi, Halforders and automorphisms of chain
    structures. \emph{J.\ Geom.} \textbf{71} (2001), 1--18.

\item %\label{}
    B.~Alinovi, H.~Karzel, Halfordered sets, halfordered chain structures
    and splittings by chains. \emph{J.\ Geom.} \textbf{75} (2002), 15--26.

\item %\label{}
    L.~Giuzzi, H.~Karzel, Co-Minkowski spaces, their reflection structure
    and K-loops. \emph{Discrete Math.} \textbf{255} (2002), 161--179.

\item %\label{}
    B.~Im, H.~Karzel, H.-J.~Ko, Webs with rotation and reflection
    properties and their relations with certain loops. \emph{Abh.\ Math.\
    Sem.\ Univ.\ Hamburg} \textbf{72} (2002), 9--20.

\item %\label{}
    H.~Karzel, S.~Pianta, E.~Zizioli, K-loops derived from Frobenius
    groups. \emph{Discrete Math.} \textbf{255} (2002), 225--234.

\item %\label{}
    H.~Karzel, S.~Pianta, E.~Zizioli, Loops, reflection structures and
    graphs with parallelism. \emph{Results Math.} \textbf{42} (2002),
    74--80.

\item %\label{}
    H.~Karzel, M.~Marchi, Relations between the K-loop and the defect of an
    absolute plane. \emph{Results Math.} \textbf{47} (2005), 305--326.

\item %\label{}
    H.~Karzel, S.~Pianta, Left loops, bipartite graphs with parallelism and
    bipartite involution sets. \emph{Abh.\ Math.\ Sem.\ Univ.\ Hamburg}
    \textbf{75} (2005), 203--214.

\item %\label{}
    H.~Karzel, S.~Pianta, E.~Zizioli, From involution sets, graphs and
    loops to loop-nearrings. In: H.~Kiechle, A.~Kreuzer, M.~J.~Thomsen
    (eds.), \emph{Nearrings and Nearfields (Proceedings of the Conference
    on Nearrings and Nearfields, Hamburg, Germany, July 27--August 3,
    2003)}, 235--252, Springer, Dordrecht 2005.

\item %\label{}
    H.~Karzel, Emanuel Sperner: Leben und Werk. \emph{Mitt.\ Math.\ Ges.\
    Hamburg} \textbf{25} (2006), 23--32.

\item %\label{}
    H.~Karzel, Emanuel Sperner: Begr\"{u}nder einer neuen Ordnungstheorie.
    \emph{Mitt.\ Math.\ Ges.\ Hamburg} \textbf{25} (2006), 33--44.

\item %\label{}
    H.~Karzel, M.~Marchi, Classification of general absolute geometries
    with Lambert-Saccheri quadrangles. \emph{Matematiche (Catania)}
    \textbf{61} (2006), 27--36.

\item %\label{}
    H.~Karzel, M.~Marchi, Vectorspacelike representation of absolute
    planes. \emph{J.\ Geom.} \textbf{86} (2006), 81--97 (2007).

\item %\label{}
    H.~Karzel, S.~Pianta, E.~Zizioli, Polar graphs and corresponding
    involution sets, loops and Steiner triple systems. \emph{Results Math.}
    \textbf{49} (2006), 149--160.

\item %\label{}
    H.~Karzel, Loops related to geometric structures. \emph{Quasigroups
    Related Systems} \textbf{15} (2007), 47--76.

\item %\label{}
    H.~Karzel, M.~Marchi, S.~Pianta, Legendre-like theorems in a general
    absolute geometry. \emph{Results Math.} \textbf{51} (2007), 61--71.

\item %\label{}
    H.~Karzel, J.~Kosiorek, A.~Matra\'{s}, Properties of auto- and
    antiautomorphisms of maximal chain structures and their relations to
    i-perspectivities. \emph{Results Math.} \textbf{50} (2007), 81--92.

\item %\label{}
    H.~Karzel, M.~Marchi, Introduction of measures for segments and angles
    in a general absolute plane. \emph{Discrete Math.} \textbf{308} (2008),
    220--230.

\item %\label{}
    H.~Karzel, S.~Pianta, Binary operations derived from symmetric
    permutation sets and applications to absolute geometry. \emph{Discrete
    Math.} \textbf{308} (2008), 415--421.

\item %\label{}
    H.~Karzel, K.~S\"{o}rensen, Lambert-Saccheri quadrangles. \emph{J.\ Geom.}
    \textbf{91} (2009), 61--62.

\item %\label{}
    H.~Karzel, J.~Kosiorek, A.~Matra\'{s}, Automorphisms of symmetric and
    double symmetric chain structures. \emph{Results Math.} \textbf{55}
    (2009), 401--416.

\item %\label{}
    H.~Karzel, J.~Kosiorek, A.~Matra\'{s}, Ordered symmetric Minkowski
    planes~I. \emph{J.\ Geom.} \textbf{93} (2009), 116--127.

\item %\label{}
    H.~Karzel, M.~Marchi, S.~Pianta, Three-reflection theorems in the
    hyperbolic plane. In: \emph{Trends in Incidence and Galois Geometries:
    A tribute to Giuseppe Tallini}, volume~19 of \emph{Quad.\ Mat.},
    127--140, Dipartimento di Matematica, Seconda Universit\`{a} di Napoli,
    Caserta 2009.

\item %\label{}
    H.~Karzel, J.~Kosiorek, A.~Matra\'{s}, Symmetric Minkowski planes
    ordered by separation. \emph{J.\ Geom.} \textbf{98} (2010), 115--125.

\item %\label{}
    H.~Karzel, M.~Marchi, S.~Pianta, The defect in an invariant reflection
    structure. \emph{J.\ Geom.} \textbf{99} (2010), 67--87.

\item %\label{}
    H.~Karzel, J.~Kosiorek, A.~Matra\'{s}, Point symmetric $2$-structures.
    \emph{Results Math.} \textbf{59} (2011), 229--237.

\item %\label{}
    H.~Karzel, M.~Marchi, S.-G.~Taherian, Elliptic reflection structures,
    K-loop derivations and triangle-inequality. \emph{Results Math.}
    \textbf{59} (2011), 163--171.

\item %\label{}
    H.~Karzel, K.~S\"{o}rensen, Metric kinematic planes and their
    representation by $\nu$-local systems. \emph{J.\ Geom.} \textbf{100}
    (2011), 85--103.

\item %\label{}
    H.~Karzel, S.-G.~Taherian, Reflection spaces, partial K-loops and
    K-loops. \emph{Results Math.} \textbf{59} (2011), 213--218.

\item %\label{}
    F.~Bonenti, H.~Karzel, M.~Marchi, Absolute planes with elliptic
    congruence. \emph{Mitt.\ Math.\ Ges.\ Hamburg} \textbf{32} (2012),
    123--143.

\item %\label{}
    H.~Karzel, S.-G.~Taherian, Reflection spaces and corresponding
    kinematic structures. \emph{Results Math.} \textbf{63} (2013),
    597--610.

\item %\label{}
    H.~Karzel, J.~Kosiorek, A.~Matra\'{s}, A representation of a point
    symmetric $2$-structure by a quasi-domain. \emph{Results Math.}
    \textbf{65} (2014), 333--346.

\item %\label{}
    H.~Karzel, J.~Kosiorek, A.~Matra\'{s}, Symmetric $2$-structures, a
    classification. \emph{Results Math.} \textbf{65} (2014), 347--359.

\item %\label{}
    H.~Karzel, S.~Pianta, S.~Pasotti, A class of fibered loops related to
    general hyperbolic planes. \emph{Aequationes Math.} \textbf{87} (2014),
    31--42.

\item %\label{}
    H.~Karzel, S.~Pianta, M.~Rostamzadeh, S.-G.~Taherian, Classification of
    general absolute planes by quasi-ends. \emph{Aequationes Math.}
    \textbf{89} (2015), 863--872.

\item %\label{}
    H.~Karzel, Loops related to reflection geometries. \emph{God.\
    Sofi\u{\i}. Univ.\ ``Sv. Kliment Okhridski''. Fak.\ Mat.\ Inform.}
    \textbf{103} (2016), 33--38.

\item %\label{}
    H.~Karzel, S.-G.~Taherian, Elliptic reflection spaces. \emph{Results
    Math.} \textbf{69} (2016), 1--10.

\item %\label{}
    H.~Karzel, S.-G.~Taherian, Groups with a ternary equivalence relation.
    \emph{Aequationes Math.} \textbf{92} (2018), 415--423.

\item %\label{}
    H.~Karzel, S.-G.~Taherian, Properties of reflection geometries and the
    corresponding group spaces. \emph{Results Math.} \textbf{74} (2019),
    Paper No. 99, 11~pp.

\item %\label{}
    H.~Karzel, Erinnerungen an Heinrich Wefelscheid. \emph{Mitt.\ Math.\
    Ges.\ Hamburg} \textbf{40} (2020), 11--14.
\end{enumerate}

\section*{Books}

\begin{enumerate}
%% local re-definition of counter
\renewcommand{\labelenumi}{{\rm B\arabic{enumi}.}}%

\item %\label{}
    H.~Karzel, K.~S\"{o}rensen, D.~Windelberg, \emph{Einf\"{u}hrung in die
    Geometrie}. Vandenhoeck \& Ruprecht, G\"{o}ttingen 1973.

\item %\label{}
    H.~Karzel, K.~S\"{o}rensen (eds.), \emph{Wandel von Begriffsbildungen in
    der Mathematik}. Wissenschaftliche Buchgesellschaft, Darmstadt 1984.

\item %\label{}
    H.~Karzel, H.-J.~Kroll, \emph{Geschichte der Geometrie seit Hilbert}.
    Wissenschaftliche Buchgesellschaft, Darmstadt 1988.

\item %label{}
    W.~Benz, H.~Karzel, A.~Kreuzer (eds.), \emph{Emanuel Sperner:
     Gesammelte Werke}. Heldermann Verlag, Lemgo 2005.

\item %\label{}
    G.~Graumann, H.~Karzel, \emph{Metric Planes. A Group Theoretical
    Foundation}. WTM-Verlag, M\"{u}nster 2018.
\end{enumerate}

\section*{}
\par\noindent
    Hans Havlicek\\
    Institut f\"{u}r Diskrete Mathematik und Geometrie\\
  	Technische Universit\"{a}t Wien\\
    Wiedner Hauptstra{\ss}e 8--10\\
    1040 Wien\\
    Austria\\
    \texttt{havlicek@geometrie.tuwien.ac.at} %%
\bigskip
\par\noindent
    Alexander Kreuzer\\
    Fachbereich Mathematik\\
    Universit\"{a}t Hamburg\\	
    Bundesstra{\ss}e 55\\
    20146 Hamburg\\
    Germany\\
    \texttt{kreuzer@math.uni-hamburg.de}
\bigskip
\par\noindent
    Hans-Joachim Kroll\\
    Department of Mathematics\\
    TU M\"{u}nchen\\
    Boltzmannstra{\ss}e 3\\
    85748 Garching bei M\"{u}nchen\\
    Germany\\
    \texttt{kroll@ma.tum.de}
\bigskip
\par\noindent
    Kay S\"{o}rensen\\
    Department of Mathematics\\
    TU M\"{u}nchen\\
    Boltzmannstra{\ss}e 3\\
    85748 Garching bei M\"{u}nchen\\
    Germany\\
    \texttt{soeren@ma.tum.de}
\end{document}